\documentclass{amsart}

\newtheorem{theorem}{Theorem}[section]
\newtheorem*{question}{Question}
\newtheorem{lemma}[theorem]{Lemma}

\theoremstyle{definition}

\theoremstyle{remark}

\numberwithin{equation}{section}

\newcommand{\ra}{\rightarrow}

\newcommand \bl {\begin{lemma}}
\newcommand \el {\end{lemma}}
\newcommand \bt {\begin{theorem}}
\newcommand \et {\end{theorem}} 
\newcommand \mc{\mathcal}

\newcommand \Ric{\mbox{Ric}}
\newcommand \Vol{\mbox{Vol}}

\newcommand \htop{h_{top}}
\newcommand \bp{\begin{proof}}
\newcommand \ep{\end{proof}}

\begin{document}
\title[Topological entropy under normalised Ricci Flow]{A criterion for topological entropy to decrease under normalised Ricci flow}
\begin{abstract}
In 2004, Manning showed that the topological entropy of the geodesic flow for a surface of negative curvature decreases as the metric evolves under the normalised Ricci flow. It is an interesting open problem, also due to Manning, to determine to what extent such behaviour persists for higher dimensional manifolds. In this short note, we describe the problem and give a strong criterion under which monotonicity of the topological entropy can be established for a short time. In particular, the criterion applies to metrics of negative sectional curvature which are in the same conformal class as a metric of constant negative sectional curvature.
\end{abstract}

\author{Daniel J. Thompson}
\address{Department of Mathematics, Pennsylvania State University, University Park, State College, PA 16802, USA}
\email{thompson@math.psu.edu}
\maketitle
\section{Introduction}
We consider the normalised Ricci flow (NRF) $\{g_t\}$ for an initial metric $g_0$ on an $n$ dimensional compact manifold $M$ which admits a metric of constant negative sectional curvature $g_C$. We let $\htop (g)$ denote the topological entropy of the geodesic flow for a metric $g$ and we consider the map
\[
t \mapsto \htop (g_t).
\]
Let $\mc S(g)$ denote the total scalar curvature functional. We show that if the initial metric $g_0$ has negative sectional curvature (or more generally, if  $g_0$ is an Anosov metric), and $\mc S(g_0) \geq \mc S(g_C)$, then there exists a small time interval $[0, \epsilon)$ on which the map $t \mapsto \htop (g_t)$ is decreasing. In particular, the condition on $\mc S(g_0)$ holds when $g_0$ and $g_C$ are in the same conformal class \cite{Schoen}. This result gives some progress towards a question of Manning \cite{Ma}, who proved that  the map $t \mapsto \htop (g_t)$ is decreasing for $t \in [0, \infty)$ for a surface of negative Gauss curvature, and asked if any analogue is true in higher dimensions.  Manning asked four other natural questions about the relationship between various dynamical quantities and Ricci flow. To date, Manning's questions remain open, although some related progress was made by Jane \cite{Jan}.


To motivate this work, we recall some dynamical and geometric interpretations of topological entropy. The quantity $\htop(g_t)$ can be interpreted as a measure of the exponential complexity of the orbit structure of the geodesic flow on $(M, g_t)$. For a precise description of this interpretation in the case of a discrete time dynamical system see \cite[page 172]{Wa}. In light of this, we can interpret the monotonicity of topological entropy along the normalised Ricci flow as a decrease in the complexity of the orbit structure of the geodesic flow. 
Furthermore, Bowen \cite{Bo4} showed that for $g$ with negative sectional curvatures,
\begin{equation} \label{def1}
\htop (g) = \lim_{t \ra \infty} \frac{1}{t} \log P(t),
\end{equation} where $P(t)$ is the number of closed geodesics of period at most $t$. 
Thus, we are describing a situation where the NRF decreases the exponential growth rate of the number of closed geodesics. For metrics with negative sectional curvature (or more generally, metrics with no conjugate points), $\htop (g)$ can be formulated as 
\begin{equation}  \label{def2}
\htop(g) = \lim_{R \ra \infty} \frac{1}{R} \log \Vol(B(x, R)),
\end{equation}
where $B(x, R)$ is the ball of radius $R$ in the universal cover of $(M,g)$, and the choice of $x$ is arbitrary \cite{Man}, \cite{FM}. This quantity is often called the volume entropy. It is interesting to compare this definition with that of the asymptotic volume ratio, where one looks for a polynomial growth rate for $\Vol(B(x, R))$. Results on the evolution of asymptotic volume ratio under Ricci flow for non-negative curvature manifolds are well known in the Ricci flow literature \cite{Ham}, \cite{MT}. 
The volume entropy can be thought of as a negative curvature analogue of the asymptotic volume ratio, so this provides further evidence that studying the evolution of $\htop(g)$ under a geometric flow is a natural thing to do. 

This note constitutes a first step in establishing a relationship between topological entropy and geometric flows for higher dimensional manifolds.  It would be desirable to understand the evolution of topological entropy under normalised Ricci flow for a larger class of metrics than those considered here,  and we discuss what one could hope to establish at the end of this note. 

\section{Results}
Let $M$ be a compact manifold which admits a metric $g_{C}$ of constant sectional curvature. We consider the family of metrics $\{g_t\}$  defined by an initial metric $g_0$ and the normalised Ricci flow (NRF) equation
\begin{equation}  \label{def3}
\frac{\partial g}{\partial t} = -2 \Ric +\frac{2}{n} \frac {\int S d \Vol}{\Vol(M)} g,
\end{equation}
where $\Ric$ is the Ricci curvature and $S$ is the scalar curvature \cite{Top}, \cite{Ham}, \cite{CK}. Of course, $\Ric$, $S$ and $d \Vol$ depend on the metric, and we write, for example, $S_g$ when we want to highlight this dependence.  We call $\{g_t\}$ the NRF for $g_0$. The second term in the NRF equation (\ref{def3}) serves to fix the volume of the metric as it evolves. We write $M_t$ for the manifold $(M, g_t)$ and $SM_t$ for the unit tangent bundle of $(M, g_t)$. Similarly, for the constant curvature metric $g_C$, we write $M_C$ and $SM_C$. We
make a standing assumption that the initial metric has the same volume as $g_C$ and is not isometric to $g_C$. We will consider initial metrics that have everywhere negative sectional curvature or are Anosov. The Anosov metrics are those for which the geodesic flow is an Anosov flow, and include those with everywhere negative curvature \cite[\S17.4]{KH}. 

We consider the map $t \mapsto \htop (g_t)$, where $\htop(g)$ is defined by (\ref{def1}), or equivalently (\ref{def2}). Our results are as follows.

\bt \label{main}
Suppose $n \geq3$ and $(M^n, g_C)$ is a compact Riemannian manifold with constant negative sectional curvature. Let $g_0$ be a metric of negative sectional curvature (or, more generally, an Anosov metric) which belongs to the conformal class of $g_C$. There exists $\epsilon$ such that $t \mapsto \htop (g_t)$ is decreasing for $t \in [0, \epsilon)$.  
\et

Let $\mc S(g_t) = \int_{M_t} S d \Vol$ denote the total scalar curvature functional. 
\bt \label{general}
Let $\{g_t\}$ be a NRF on a compact Riemannian manifold which admits a metric $g_C$ of constant negative sectional curvature. For a fixed time $t$, suppose $g_t$ has everywhere negative sectional curvature (or is an Anosov metric), and that $-\mc S(g_t) \leq - \mc S (g_C)$, then
\[
\frac{d}{ds} (\htop (g_s))  |_{s=t} < 0. 
\]
\et
Theorem \ref{main} is an immediate consequence of theorem \ref{general}. This is because
\[
\mc S (g_C) = \inf \{ \mc S (g) : g \mbox{ is in the conformal class of } g_C, \Vol(g) =\Vol(g_C) \},
\]
see \cite[proposition 1.4]{Schoen}. Thus, the hypothesis of theorem \ref{general} is satisfied at the initial metric in theorem \ref{main}.

Our proof of theorem \ref{general} follows the same strategy that Manning used in the case of surfaces, and is based on the Katok-Kneiper-Weiss formula \cite{KKW} for the derivative of topological entropy along a path of metrics  (the formula applies in all dimensions). 
We require that the topological entropy is mimimised at $g_C$. For $n \geq3$, this is the celebrated Besson-Courtois-Gallot theorem \cite{BCG}. The result was originally conjectured by Katok, who proved it for surfaces \cite{Ka} and in the conformal class of $g_C$. The other key ingredient is a theorem of Freire and Ma\~n\'e \cite{FM} which relates  $\int_{SM} \Ric (v,v) d \mu(v)$ and $h_{\mu}$ for any probability measure $\mu$ on $SM$ which is invariant under the geodesic flow (this is a corollary of the main result of \cite{FM} which gives a relationship between certain solutions of the matrix Riccati equation and Lyapanov exponents of invariant measures).  
\bp

The Katok-Kneiper-Weiss formula \cite{KKW} says that 
\[
\frac{d}{ds} (\htop (g_s))  |_{s=t} = -\frac{\htop(g_t)}{2} \int_{SM_{t}} (\partial / \partial s) |_{s=t} g_s (v, v) d \mu_t (v),
\]
where  $t \mapsto g_t$ is a $C^2$ path of Anosov Riemannian metrics. 
Here, $\mu_t$ is the Bowen-Margulis measure for $M_t$, which is the unique probability measure that satisfies $h_{\mu_t} = \htop (g_t)$, where $h_{\mu_t}$ is the measure theoretic entropy of $\mu_t$ with respect to the geodesic flow on $M_t$. On substituting the NRF equation into the Katok-Kneiper-Weiss formula, the problem reduces to showing that
\[
\int_{SM_{t}}  -\Ric d \mu_t - \frac{1}{n \Vol(M)} \int_{M_t} -S d \mbox{Vol} > 0.
\]

Corollary II.1 of Freire and Ma\~n\'e \cite{FM} states that for a metric without conjugate points on $M$  and any probability measure $\nu$ on $SM$ that is invariant under the geodesic flow, then
\[
(n-1) \int_{SM} -\Ric d \nu \geq h_\nu^2,
\]
We note that in \cite{FM}, the authors define $\Ric$ as the usual Ricci curvature divided by $(n-1)$, which means that their statement of corollary II.1 appears to differ from the inequality above by a factor of $(n-1)$. Applying this to the Bowen-Margulis measure $\mu_t$, we have 
\[
(n-1) \int_{SM_t} -\Ric d \mu_t \geq h_{\mu_t}^2 > h_{\mu_C}^2,
\]
where the second inequality is an application of the Besson-Courtois-Gallot theorem \cite{BCG} and $\mu_C$ is the Bowen-Margulis measure (equivalently Liouville measure) for $(M, g_C)$. Thus, 
\[
\int_{SM} -\Ric d \mu_t - \frac{1}{n-1} h_{\mu_C}^2 > 0.
\]
At $g_C$, suppose the constant negative curvature is $-K$. Then $h_{\mu_C}^2 = (n-1)^2 K$ and $S_{g_C} = -n(n-1) K$. Therefore,
\[
h_{\mu_C}^2 = \frac{(n-1)}{n \Vol(M)} \int_{M_C} - S_{g_C} d \Vol.
\]
We have assumed that $-\mc S(g_t) \leq - \mc S(g_C)$ and thus

\[
\frac{1}{n \Vol(M)}  \int_{M_t} -S_{g_t} d \mbox{Vol}\leq \frac{1}{n \Vol(M)}  \int_{M_C} -S_{g_C} d \mbox{Vol} =  \frac{1}{n-1}h_{\mu_C}^2.
\]
Therefore, 
\[
\int_{SM_{t}}  -\Ric d \mu_t - \frac{1}{n \Vol(M)} \int_{M_t} -S d \mbox{Vol} \geq \int_{SM_t} -\Ric d \mu_t - \frac{1}{n-1} h_{\mu_C}^2 > 0,
\]
and thus$ \frac{d}{ds} (\htop (g_s))  |_{s=t} < 0.$ 
\ep

\section{Remarks}
\subsection{Future directions} We discuss what more one might hope to say in higher dimensions.
\begin{question} [Manning]
Does the map $t \mapsto \htop (g_t)$ decrease for $t\in[0, \infty)$ when $\{g_t\}$ is a NRF with initial metric in a neighbourhood of $g_C$?
\end{question}
The question certainly makes sense because a result of Ye \cite{Ye} states that if $g_0$ is in a $C^2$ neighbourhood of $g_C$, then the NRF of $g_0$ exists for all time, preserves the negativity of the sectional curvature, and converges to $g_C$, which by Besson-Courtois-Gallot, is the metric which minimises topological entropy. It seems that this question is rather a subtle problem. Theorem \ref{general} does not apply because although $g_C$ minimises the total scalar curvature functional in the conformal class, we expect that $\mc S(g) < \mc S(g_C)$ for a metric $g$ given by a perturbation of $g_C$ in a transversal direction \cite{Schoen}. 
Further progress will require more sophisticated techniques than the short argument presented here.
\subsection{Evolution of topological entropy under the unnormalised Ricci flow} \label{URF} We first point out that the question of how $\htop(g_t)$ behaves when $g_t$ evolves by the Ricci flow (RF) equation $\frac{\partial g}{\partial t} = -2 \Ric $ is not a good one. This is because $\htop(g_t)$ is not scale invariant. For example, it decreases under a homothetic expansion of the manifold. Since the Ricci flow increases volume when the curvature is everywhere negative, we would expect that it is easier to show that $\htop(g_t)$ is decreasing under RF than under NRF. This is the case, because as long as $\Ric$ is negative, it follows directly from the Katok-Kneiper-Weiss formula that
\[
\frac{d}{ds} (\htop (g_s))  |_{s=t} = \htop(g_t) \int_{SM_{t}} \Ric (v,v) d \mu_t (v) < 0.
\]
If the initial metric is in a neighbourhood of $g_C$, then by Ye \cite{Ye}, the flow exists for all time, and $\Ric_{g_t}(v,v)$ stays negative. Thus, $t \mapsto \htop (g_t)$ decreases for $t\in[0, \infty)$.  As discussed in the introduction, we can interpret the monotonicity of topological entropy along the {\bf normalised} Ricci flow as a decrease in the complexity of the orbit structure of the geodesic flow. However, a priori, monotonicity of topological entropy along the {\bf unnormalised} Ricci flow could just be due to the expansion of the volume along the flow, and would not necessarily reveal anything about the complexity of the orbit structure.
\section*{Acknowledgements}
I would like to thank Xiaodong Wang for making me aware of an error in an earlier version of this work, which has now been resolved.

\bibliographystyle{plain}
\bibliography{master}
\end{document}